\documentclass{article}
\usepackage{amsmath,amssymb,latexsym,braket,hyperref,cleveref,url}

\newcommand{\merely}[1]{\ensuremath{\Vert#1\Vert}}
\newcommand{\trunc}[2]{\ensuremath{\Vert#2\Vert_{#1}}}
\newcommand{\tr}[1]{\ensuremath{|#1|}}
\newcommand{\U}{\ensuremath{\mathcal{U}}}
\newcommand{\ct}{%
  \mathchoice{\mathbin{\raisebox{0.5ex}{$\displaystyle\centerdot$}}}%
             {\mathbin{\raisebox{0.5ex}{$\centerdot$}}}%
             {\mathbin{\raisebox{0.25ex}{$\scriptstyle\,\centerdot\,$}}}%
             {\mathbin{\raisebox{0.1ex}{$\scriptscriptstyle\,\centerdot\,$}}}
}
\newcommand{\mink}{\mathsf{Mink}}
\usepackage{fullpage}
\usepackage{xspace,microtype}

\newcommand{\hottonly}{\textsc{h}o\textsc{tt}}
\newcommand{\ufonly}{\textsc{uf}}
\newcommand{\hott}{\hottonly/\ufonly}
\newcommand{\zfc}{\textsc{zfc}}
\newcommand{\etcs}{\textsc{etcs}}
\newcommand{\hit}{\textsc{hit}}
\newcommand{\hits}{\textsc{hit}s}

\crefformat{section}{\S#2#1#3}
\Crefformat{section}{Section~#2#1#3}
\crefrangeformat{section}{\S\S#3#1#4--#5#2#6}
\Crefrangeformat{section}{Sections~#3#1#4--#5#2#6}
\crefmultiformat{section}{\S\S#2#1#3}{ and~#2#1#3}{, #2#1#3}{ and~#2#1#3}
\Crefmultiformat{section}{Sections~#2#1#3}{ and~#2#1#3}{, #2#1#3}{ and~#2#1#3}
\crefrangemultiformat{section}{\S\S#3#1#4--#5#2#6}{ and~#3#1#4--#5#2#6}{, #3#1#4--#5#2#6}{ and~#3#1#4--#5#2#6}
\Crefrangemultiformat{section}{Sections~#3#1#4--#5#2#6}{ and~#3#1#4--#5#2#6}{, #3#1#4--#5#2#6}{ and~#3#1#4--#5#2#6}

\let\ea\expandafter
\def\alwaysmath#1{\ea\ea\ea\global\ea\ea\ea\let\ea\ea\csname your@#1\endcsname\csname #1\endcsname
  \ea\def\csname #1\endcsname{\ensuremath{\csname your@#1\endcsname}\xspace}}
\alwaysmath{infty}

\title{Homotopy Type Theory:\\A synthetic approach to higher equalities}
\author{Michael Shulman\thanks{This material is based on research sponsored by The United States Air Force Research Laboratory under agreement number FA9550-15-1-0053.  The U.S. Government is authorized to reproduce and distribute reprints for Governmental purposes notwithstanding any copyright notation thereon.  The views and conclusions contained herein are those of the author and should not be interpreted as necessarily representing the official policies or endorsements, either expressed or implied, of the United States Air Force Research Laboratory, the U.S. Government, or Carnegie Mellon University.}}
\date{}

\begin{document}
\maketitle

\tableofcontents

\section{Introduction}
\label{sec:introduction}

Ask an average mathematician or philosopher today about the foundations of mathematics, and you are likely to receive an answer involving set theory: an apparent consensus in marked contrast to the foundational debates of the early 20th century.
Now, at the turn of the 21st century, a new theory has emerged to challenge the foundational ascendancy of sets.
Arising from a surprising synthesis of constructive intensional type theory and abstract homotopy theory, Homotopy Type Theory and Univalent Foundations (\hott) purports to represent more faithfully the everyday practice of mathematics, but also provides powerful new tools and a new paradigm.
So far, its concrete influence has been small, but its potential implications for mathematics and philosophy are profound.

There are many different aspects to \hott,\footnote{Though \hottonly\ and \ufonly\ are not identical, the researchers working on both form a single community, and the boundary between them is fluid.  Thus, I will not attempt to distinguish between them, even if it results in some technically incorrect statements.} % that ought really to apply only to one or the other.}
but in this chapter I will focus on its use as a foundation for mathematics.
% mainly on one: it offers a new foundation for mathematics, which could be used instead of set theory. % \footnote{Although as we will see, in a certain sense it is a strict \emph{generalization} of set theory.}
Like set theory, it proposes to found mathematics on a notion of \emph{collection}, but its collections (called \emph{types}) behave somewhat differently.
% For instance, rather than being pre-existing objects that are merely grouped into a set, the elements of a type have no existence apart from that type.
The most important difference is that in addition to having elements as sets do, the types of \hott\ come with further collections of \emph{identifications} between these elements (i.e.\ ways or reasons that they are equal).
These identifications form a structure that modern mathematicians call an \emph{$\infty$-groupoid} or \emph{homotopy type}, which is a basic object of study in homotopy theory and higher category theory; thus \hott\ offers mathematicians a new approach to the latter subjects.

Of greater importance philosophically, however, is \hott's proposal that such types can be the fundamental objects out of which mathematics and logic are built.
In other words, \hott\ suggests that whenever we mentally form a collection of things, we must \emph{simultaneously} entertain a notion of what it means for two of those things to be the same (in contrast to the position of \zfc\ that all things have an identity criterion \emph{prior} to their being collected into a set).
As stated, this is closely related to the conception of ``set'' promulgated by Bishop; but \hott\ generalizes it by allowing two things to ``be the same'' in \emph{more than one way}.
This is perhaps not a common everyday occurrence, but it is a fundamental part of category theory and thus an integral part of mathematics, including many modern theories of physics.
Thus, like other initially unintuitive ideas such as relativistic time dilation and quantum entanglement, it can be argued to be basic to the nature of reality.
The innovation of \hott\ is that this idea can be made basic to the foundational logical structure of mathematics as well, and that doing so actually \emph{simplifies} the theory.

In this chapter I will attempt to convey some of the flavor and advantages of \hott;
we will see that in addition to {expanding} the discourse of mathematics, it also represents certain aspects of \emph{current} mathematical practice more faithfully than set theory does.
In \crefrange{sec:infty-groupoids}{sec:synthetic} I will describe \hott\ very informally; in \crefrange{sec:type-theory}{sec:high-induct-types} I will discuss some of its features in a bit more detail; % about some of its features such as logic, univalence, and \hits;
and in \cref{sec:erlangen} I will attempt to pull together all the threads with an example.
For space reasons, I will not be very precise, nor will I discuss the history of the subject in any depth; for more details see~\cite{hottbook}.
Other recent survey articles on \hott\ include~\cite{apw:vvu-hott,awodey:tt-and-htpy,pw:hottvvuf}.

For helpful conversations and feedback, I would like to thank (in random order) Emily Riehl, David Corfield, Dimitris Tsementzis, James Ladyman, Richard Williamson, Mart\'\i{}n Escard\'o, Andrei Rodin, Urs Schreiber, John Baez, and Steve Awodey; as well as numerous other contributors at the $n$-Category Caf\'{e} and the \hottonly\ email list, and the referees.

\section{$\infty$-groupoids}
\label{sec:infty-groupoids}

The word ``$\infty$-groupoid'' looks complicated, but the underlying idea is extremely simple, arising naturally from a careful consideration of what it means for two things to be ``the same''.
% We begin with the observation that in mathematics, and also in everyday language, we often have occasion to consider collections of things.
% The prevailing language used to describe such collections is \emph{set theory}, which uses the word ``set'' for such a collection, and supplies a rich toolbox for manipulating them.
% One writes $\{a,b,c\}$ for the set whose elements are $a$, $b$, and $c$.
Specifically, it happens frequently in mathematics that we want to define a collection of objects that are determined by some kind of ``presentation'', but where ``the same'' object may have more than one presentation.
As a simple example, % consider the problem of defining real numbers.
if we try to define a \emph{real number} to be an infinite decimal expansion\footnote{Like any mathematical object, there are many equivalent ways to define the real numbers.
This specific definition is rarely used in mathematics for technical reasons, but it serves as a good illustration, and the common definition of real numbers using Cauchy sequences has exactly the same issues.}
such as $\pi = 3.14159\cdots$, we encounter the problem that (for instance)
\[ 0.5=0.50000\cdots \qquad\text{and}\qquad 0.4\overline{9}=0.49999\cdots \]
are distinct decimal expansions but ought to {represent the same real number}.
Therefore, ``the collection of infinite decimal expansions'' is not a correct way to define ``the collection of real numbers''.

If by ``collection'' we mean ``set'' in the sense of \zfc, then we can handle this by defining a real number to be a \emph{set} of decimal expansions that all ``define the same number'', and which is ``maximal'' in that there are no \emph{other} expansions that define the same number.
Thus, one such set is $\{0.5, 0.4\overline{9}\}$, and another is $\{0.\overline{3}\}$.
% Many other such sets are singletons, such as $\{0.\overline{3}\}$, since duplication is only possible for expansions ending in all $0$'s or all $9$'s.
% Thus, the set of real numbers might look like this:
% \[ \Big\{ \; \{0.5, 0.4\overline{9}\} ,\; \{0.\overline{3}\} ,\; \dots \; \Big\} \]
% The sets such as $\{0.5, 0.4\overline{9}\}$ and $\{0.\overline{3}\}$ 
These sets are \emph{equivalence classes}, and the information about which expansions define the same number is an \emph{equivalence relation} % on the collection of decimal expansions.
(a binary relation $\sim$ such that $x\sim x$, if $x\sim y$ then $y\sim x$, and if $x\sim y$ and $y\sim z$ then $x\sim z$).
% The equivalence classes are then the sets of the form $\{ y \mid y\sim x \}$ for some $x$.)
The set of equivalence classes is the \emph{quotient} of the equivalence relation.

Similarly, Frege~\cite[\S68]{frege:grundlagen-arith} defined the \emph{cardinality} of a set $X$ to be (roughly, in modern language) the set of all sets that are related to $X$ by a bijection.
%, i.e.\ related to it by a one-to-one correspondence (bijection).
Thus for instance $0$ is the set of all sets with no elements, $1$ is the set of all singleton sets, and so on.
These are exactly the equivalence classes for the equivalence relation of bijectiveness. %\footnote{There are some issues with Frege's definition that I will not address here, such as the fact that in \zfc, all of his nonzero numbers are proper classes.}
That is, we consider a cardinal number to be ``presented'' by a set having that cardinality, with two sets presenting the same cardinal number just when they are bijective.

An example outside of pure mathematics involves Einstein's theory of general relativity, %(GR)
in which the universe is represented by a differentiable manifold % (a mathematical object)
with a metric structure. % (representing gravity) satisfying a dynamical equation.
% A cornerstone of this theory is \emph{general covariance}, the principle that physical reality does not depend on coordinate systems.\footnote{To be precise, although general covariance was historically important to Einstein's discovery of GR, Kretschmann~\cite{kretschmann:covariance} pointed out that
%   % with sufficient mathematical ingenuity
%   nearly \emph{any} physical theory can be formulated covariantly, and nowadays general covariance is usually not regarded as unique to GR; see~\cite{norton:covariance} for a historical survey.
%   So the point I am making here is not really specific to GR, but it is perhaps easiest to understand in that case.% , since many people are more used to thinking of GR covariantly than, say, Newtonian theory.
% }
% In particular, % this leads to the conclusion that
% ``{General covariance}'' implies that
In this theory, if two manifolds are \emph{isomorphic} respecting their metric structure, then they represent the same physical reality.
(An isomorphism of manifolds is often called a ``diffeomorphism'', and if it respects the metric it is called an ``isometry''.)
% , in the same way that $0.5$ and $0.4\overline{9}$ ought to represent the same real number.
Thus we find for instance in~\cite[\S1.3]{sachs-wu:relativity} that
\begin{quote}
  \small
  A general relativistic \emph{gravitational field} $[(M,\mathbf{g})]$ is an equivalence class of spacetimes [manifolds $M$ with metrics $\mathbf{g}$] where the equivalence is defined by \dots\ 
  % orientation and time-orientation-preserving
  isometries.
\end{quote}
This sort of situation, where multiple mathematical objects represent the same physical reality, is common in modern physics, and the mathematical objects (here, the manifolds) are often called \emph{gauges}.%
\footnote{Whether or not general relativity should be technically considered a ``gauge theory'' is a matter of some debate, but all that matters for us is that it exhibits the same general phenomenon of multiple models.}

Definitions by equivalence classes are thus very common in mathematics and its applications, but they are not the only game in town.
A different approach to the problem of ``presentations'' was proposed by Bishop~\cite[\S1.1]{bishop:fca}:
\begin{quote}
  \small
  A set is defined by describing exactly what must be done in order to construct an element of the set and what must be done in order to show that two elements are equal.
\end{quote}
In other words, according to Bishop, a \emph{set} is a collection of things \emph{together with} the information of when two of those things are equal (which must be an equivalence relation).\footnote{Although Bishop's goal was to give a constructive treatment of mathematics, this notion of ``set'' is meaningful independently of whether one's logic is constructive or classical.}
Thus the real numbers would \emph{be} infinite decimal expansions, but ``the set of real numbers'' would include the information that (for instance) $0.5$ and $0.4\overline{9}$ are the same real number.
% Similarly, a gravitational field would \emph{be} a manifold, but the set of gravitational fields would come with the information that isometric manifolds define the same gravitational field.
One advantage of this is that if we are given ``a real number'', we never need to worry about \emph{choosing} a decimal expansion to represent it.
(Of course, for decimal expansions there are canonical ways to make such a choice, but in other examples there are not.)

% These examples show the broad application of definitions by equivalence classes in mathematics.
% For the most part, such definitions work, but they have a few drawbacks.
% One is that it looks a bit artificial to define a real number to be a \emph{set} of decimal expansions rather than a single one.
% A more serious problem is that
% if we are given ``a real number'' % $x\in \mathbb{R}$
% or ``a gravitational field'' and we want to do any computations with it, we need to choose an element of it (called a \emph{representative} of the equivalence class).
% For decimal expansions, this is not difficult, but in the other two examples above the equivalence classes are infinite, so it may not be obvious how to do this. % choose a particular element (and may even require even the infamous \emph{axiom of choice}).

As a much older example of this style of definition,
%While Bishop may have been the first to formulate this conception of ``set'' explicitly, the idea is much older.
%For instance,
in Euclid's \emph{Elements} we find:
\begin{quotation}
  \noindent\small
\textbf{Definition 4.} Magnitudes are said to \emph{have a ratio} to one another which can, when multiplied, exceed one another.\\
\textbf{Definition 5.} Magnitudes are said to be \emph{in the same ratio}, the first to the second and the third to the fourth, when, if any equimultiples whatever are taken of the first and third, and any equimultiples whatever of the second and fourth, the former equimultiples alike exceed, are alike equal to, or alike fall short of, the latter equimultiples respectively taken in corresponding order.
\end{quotation}
That is, Euclid first defined how to \emph{construct} a ratio,  % (give a pair of magnitudes, some multiple of each of which is greater than the other),
and then secondly he defined when two ratios are \emph{equal}, exactly as Bishop says he ought.

% One can even argue that Bishop's proposal is closer to the way we categorize objects outside of mathematics.
% For instance, the division of animals into species such as \textit{Canis familiaris} is not defined by enumerating all the animals belonging to each species (the equivalence classes), but by specifying what it means for two given animals to belong to the same species (e.g.\ they can interbreed and produce fertile offspring).
% Less formally, psychological \emph{exemplar theory} suggests that mentally, a conceptual category such as ``dog'' is learned by way of examples, including information about how instances of the same concept may differ from each other; not too different from specifying the elements of a set and when they are the same.

On its own, Bishop's conception of set is not a very radical change.
%; it amounts to keeping an equivalence relation around rather than passing to the set of its equivalence classes.
But it paves the way for our crucial next step, which is to recognize that frequently there may be more than one ``reason'' why two ``presentations'' define the same object.
For example, there are two bijections between $\{a,b\}$ and $\{c,d\}$: one that sends $a$ to $c$ and $b$ to $d$, and another that sends $a$ to $d$ and $b$ to $c$.
Likewise, a pair of manifolds may be isometric in more than one way. %, and a manifold may even be isometric to itself in many ways.

This should not be confused with the question of whether there is more than one \emph{proof} that two things are the same.
Rather, the question is whether substituting one for the other in a mathematical statement or construction can yield multiple inequivalent results.
For instance, there is a predicate $P$ on $\{a,b\}$ such that $P(a)$ is true and $P(b)$ is false.
% if we regard $0\cdot 0=0$ as a (true) statement about (one element of) the set $\{0,2\}$, then transporting it along a bijection from $\{0,2\}$ to $\{1,3\}$ might yield either the true statement $1\cdot 1=1$ or the false statement $3\cdot 3=3$.
We can ``transport'' $P$ along a bijection from $\{a,b\}$ to $\{c,d\}$ to obtain a predicate $Q$ on $\{c,d\}$; but the resulting $Q$ will depend on which bijection we use.
If we use the bijection that sends $a$ to $c$ and $b$ to $d$, then $Q(c)$ will be true and $Q(d)$ will be false; but if we use the other bijection, then $Q(c)$ will be false and $Q(d)$ will be true.
Thus, $\{a,b\}$ and $\{c,d\}$ ``are the same'' in more than one way.

If a predicate or construction is left literally unchanged by this sort of substitution, it is called \emph{invariant}.
Thus, physicists speak of \emph{gauge invariance} when talking about theories with multiple mathematical models of the same reality.
More generally, a construction that ``varies appropriately'' under such substitutions (but in a way potentially dependent on the ``reason'' for sameness, as above), is called \emph{covariant}.
In particular, general relativity is said to be \emph{generally covariant}, meaning that a mathematical model of reality can be replaced by any isometric one --- but in a way dependent on the particular isometry chosen.

This behavior lies at the root of Einstein's famous \emph{hole argument}, which can be explained most clearly as follows.
Suppose $M$ and $N$ are manifolds with spacetime metrics $\mathbf{g}$ and $\mathbf{h}$ respectively, and $\phi$ is an isometry between them.
Then any point $x\in M$ corresponds to a unique point $\phi(x)\in N$, both of which represent the same ``event'' in spacetime.
Since $\phi$ is an isometry, the gravitational field around $x$ in $M$ is identical to that around $\phi(x)$ in $N$.
However, if $\psi$ is a \emph{different} isomorphism from $M$ to $N$ which does \emph{not} respect the metrics, then the gravitational field around $x$ in $M$ may be quite different from that around $\psi(x)$ in $N$.

So far, this should seem fairly obvious.
But Einstein originally considered only the special case where $M$ and $N$ happened to be the same manifold (though not with the same metric), where $\psi$ was the identity map $\mathrm{id}_M$, and where $\phi$ was the identity outside of a small ``hole''.
In this case, it seemed wrong that two metrics could be the same outside the hole but different inside of it.
The solution is clear from the more general situation in the previous paragraph:
the fact that the two metrics ``represent the same reality'' is witnessed by the isomorphism $\phi$, not $\psi$.
Thus, even for a point $x$ inside the hole, we should be comparing $\mathbf{g}$ at $x$ with $\mathbf{h}$ at $\phi(x)$, not with $\mathbf{h}$ at $\mathrm{id}_M(x)=x$.\footnote{While this description in modern language makes it clear why there is no paradox, it does obscure the reasons why for many years people \emph{thought} there was a paradox!
I will return to this in \cref{sec:erlangen}.}

This and other examples show that it is often essential to \emph{remember which} isomorphism we are using to treat two objects as the same.
The set-theoretic notion of equivalence classes is unable to do this, but Bishop's approach can be generalized to handle it.
Indeed, such a generalization is arguably already latent in Bishop's constructive phrasing: both the construction of elements and the proofs of equality are described in terms of \emph{what must be done}, so it seems evident that just as there may be more than one way to construct an element of a set, there may be more than one way to show that two elements are equal.
Bishop made no use of this possibility, but \hott\ takes it seriously.
The laws of an equivalence relation then become algebraic structure on these ``reasons for equality'': given a way in which $x=y$ and a way in which $y=z$, we must have an induced way in which $x=z$, and so on, satisfying natural axioms.
The resulting structure is called a \emph{groupoid}.
Thus, for instance, spacetime manifolds form a groupoid, in which the ways that $M=N$ are the isometries from $M$ to $N$ (if any exist).

If it should happen that for every $x$ and $y$ in some groupoid, there is \emph{at most one} reason why $x=y$, then our groupoid is essentially just a set in Bishop's sense; thus the universe of sets is properly included in that of groupoids.
This is what happens with decimal expansions: there is only one way in which $0.5$ and $0.4\overline{9}$ represent the same real number
% \footnote{The situation with cardinalities is a bit more subtle; I will return to it in \cref{sec:cardinality}.}
%This should not be confused with the question of whether there is more than one \emph{proof} that $0.5 = 0.4\overline{9}$; %, which is irrelevant to the matter at hand.
(i.e.\ in any statement or construction involving $0.5$, there is only one way to replace $0.5$ by $0.4\overline{9}$).
This is in contrast to the situation with manifolds, where using a different isomorphism $\phi$ or $\psi$ from $M$ to $N$ can result in different statements, e.g.\ one which speaks about $\phi(x)\in N$ and another about $\psi(x)\in N$.

The final step of generalization is to notice that we introduced {sets} (and generalized them to groupoids) to formalize the idea of ``collection''; but we have now introduced, for each pair of things $x$ and $y$ in a groupoid, an \emph{additional} collection, namely the ways in which $x$ and $y$ are equal. % \footnote{Lest there be any confusion, I should emphasize that this collection will often be empty.}
Thus, it seems natural that this collection should itself be a set, or more generally a groupoid; so that two ways in which $x=y$ could themselves be equal or not, and perhaps in more than one way.
Taken to its logical conclusion, this observation demands an infinite tower consisting of elements, ways in which they are equal, ways in which those are equal, ways in which \emph{those} are equal, and so on.
Together with all the necessary operations that generalize the laws of an equivalence relation, this structure is what we call an \emph{$\infty$-groupoid}.

This notion may seem very abstruse, but over the past few decades $\infty$-groupoids have risen to a central role in mathematics and even physics, starting from algebraic topology and metastasizing outwards into commutative algebra, algebraic geometry, differential geometry, gauge field theory, computer science, logic, and even combinatorics.
It turns out to be very common that two things can be equal in more than one way.

\section{Foundations for mathematics}
\label{sec:synthetic}

In \cref{sec:infty-groupoids} I introduced the notion of $\infty$-groupoid informally. %; but how are we to work with them precisely in mathematics?
At this point a modern mathematician would probably try to give a \emph{definition} of $\infty$-groupoid, such as ``an $\infty$-groupoid consists of a collection of elements, together with for any two elements $x,y$ a collection of ways in which $x=y$, and for any two such ways $f,g$ a collection of ways in which $f=g$, and so on, plus operations \dots''
Clearly, any such definition must refer to a \emph{prior} notion of ``collection'', which a modern mathematician would probably interpret as ``set''.
Such definitions of $\infty$-groupoids are commonly used, although they are quite combinatorially complicated.

However, in \cref{sec:infty-groupoids} we considered $\infty$-groupoids not as \emph{defined in terms of} sets, but as \emph{substitutes} or rather \emph{generalizations} of them.
Thus, we should instead seek a theory at roughly the same ontological level as \zfc, whose basic objects are $\infty$-groupoids.
% (Bishop did not go this far; his ``sets'' are not basic objects, but are defined in terms of a more basic notion of ``pre-set''.) % (i.e.\ a set is a pre-set equipped with an equivalence relation called ``equality'').
This is exactly what \hott\ is: a \emph{synthetic theory of $\infty$-groupoids}.\footnote{Since $\infty$-groupoids are a formalization of the idea that things can be equal in more than one way, that these ways can themselves be equal in more than one way, and so on, we may equivalently (but more informally) call \hott\ a \emph{synthetic theory of higher equalities}, as in the chapter title.}

The word ``synthetic'' here is, as usual, used in opposition to ``analytic''.
In modern mathematics, an analytic theory is one whose basic objects are defined in some other theory, whereas a synthetic theory is one whose basic objects are undefined terms given meaning by rules and axioms.
For example, \emph{analytic geometry} defines points and lines in terms of numbers; whereas \emph{synthetic geometry} is like Euclid's with ``point'' and ``line'' essentially undefined.\footnote{Euclid's \emph{Elements} as they have come down to us do contain ``definitions'' of ``point'' and ``line'', but these are not definitions in a modern mathematical sense, and more modern versions of Euclidean geometry such as Hilbert's~\cite{hilbert:geometry} do leave these words undefined.}
% Hilbert~\cite{hilbert:zahlbegriff} used the words "genetic" and "axiomatic".

% At one level, modern mathematics is characterized by a rich interplay between analytic and synthetic --- although most mathematicians would speak instead of \emph{definitions} and \emph{examples}.
% For instance, a modern geometer might define ``a geometry'' to satisfy axioms like Euclid's, and then work synthetically with those axioms; but she would also construct examples of such ``geometries'' analytically, such as with ordered pairs of real numbers.
% % Nowadays the same could be said of any other mathematical notion: groups, rings, fields, topological spaces, manifolds, algebraic varieties, etc.
% One of the first such geometers was Hilbert~\cite{hilbert:geometry}, % (see~\cite{hartshorne:geometry} for a nice modern treatment),
% who emphasized in particular that constructing an analytic example (or \emph{model}) proves the \emph{consistency} of the synthetic theory.

Thus, our first step to understanding \hott\ is that it is an axiomatic system in which ``$\infty$-groupoid'' is essentially an undefined term.
One advantage of this can already be appreciated: it allows us to say simply that for any two elements $x$ and $y$ of an $\infty$-groupoid, the ``ways in which $x=y$'' form another $\infty$-groupoid, so that $\infty$-groupoids are really the only notion of ``collection'' that we need consider.
As part of a \emph{definition} of $\infty$-groupoid, this would appear circular; %\footnote{or at least coinductive};
but as an \emph{axiom}, it is unobjectionable.
% This single change already simplifies our lives significantly: no longer do we need to think of an $\infty$-groupoid as an infinite tower of equalities and meta-equalities.

So far, this description of \hott\ could also be applied (with different terminology) to the field of mathematics called ``abstract homotopy theory''.
However, although \hott\ is strongly influenced by homotopy theory, there is more to it: as suggested above, its $\infty$-groupoids can substitute for sets as a foundation for mathematics.

% For clarity, let me state that when
% At this point I should clarify what I mean by ``foundations of mathematics''. %, since this has been used by different people to mean many different things.
When I say that a synthetic theory can be a \emph{foundation for mathematics}, I mean simply that we can encode the rest of mathematics into it somehow.\footnote{Or into some natural variant or extension of it, such as by making the logic intuitionistic or adding stronger axioms.}
This definition of ``foundation'' is reasonably precise and objective, and agrees with its common usage by most mathematicians.
A computer scientist might describe such a theory as ``mathematics-complete'', by analogy with Turing-complete programming languages (that can simulate all other languages) and NP-complete problems (that can solve all other NP problems).
For example, it is commonly accepted that \zfc\ set theory has this property.
On the other hand, category theory in its role as an organizing principle for mathematics, though of undoubted philosophical interest, is not foundational in this sense (although a synthetic form of category theory like~\cite{lawvere:catofcats} could be).

% A choice of foundational theory may of course be relevant to many other questions in the philosophy of mathematics, but this is not a precondition for its being foundational.
%
% The question of whether it is what we \emph{should} mean by a ``foundation for mathematics'' is beyond the scope of this chapter.
% If it were decided that ``foundation'' should mean something different, then we would simply have to find another word to describe this sort of theory.
%
% Of course, there are many other \emph{questions} of interest in the philosophy of mathematics, for which a choice of foundational theory may be relevant or necessary.
%
% such as what mathematical entities and structures really are, how different areas of mathematics are related, and so on.
% A foundational theory may have interesting things to say about many of these questions, and it may even be necessary to settle on some foundation before they can be asked.
%
% But it is unreasonable to \emph{require} that a theory address these questions to be considered foundational: the foundation of a building need not also be a blueprint, an architectural style, or a rooftop garden.

In particular, a synthetic theory cannot fail to be foundational because some analytic theory describes similar objects.
The fact that we \emph{can} define and study $\infty$-groupoids inside of set theory says nothing about whether a \emph{synthetic} theory of $\infty$-groupoids can be foundational.
To the contrary, in fact, it is highly \emph{desirable} of a new foundational theory that we can translate back and forth to previously existing foundations; among other things it ensures the relative consistency of the new theory.
Similarly, we cannot dismiss a new synthetic foundational theory by claiming that it ``requires some pre-existing notions'': the simple fact of being synthetic {means} that it does not.
Of course, humans always try first to \emph{understand} new ideas in terms of old ones, but that doesn't make the new ideas \emph{intrinsically} dependent on the old.
A student may learn that dinosaurs are like ``big lizards'', but that doesn't make lizards logically, historically, or genetically prior to dinosaurs.

In addition, we should beware of judging a theory to be more intuitive or fundamental merely because we are familiar with it: intuition is not fixed, but can be (and is) trained and developed.
At present, most mathematicians think of $\infty$-groupoids in terms of sets because they learned about sets early in their mathematical education; but even in its short existence the \hott\ community has already observed that graduate students who are ``brought up'' thinking in \hott\ form a direct understanding and intuition for it that sometimes outstrips that of those who ``came to it late''.
Moreover, the \zfc-like intuitions about set theory now possessed by most mathematicians and philosophers also had to be developed over time:
Lawvere~\cite{lawvere:cohesive-cantor} has pointed out that Cantor's original ``sets'' seem more like those of Lawvere's alternative set theory \etcs\ (see~\cite{lawvere:etcs-long} and McLarty's chapter in the present volume).

The point being made, therefore, is that \hott, the synthetic theory of $\infty$-groupoids, can be a foundation for mathematics in this sense.
There is quite an easy proof of this: we have already seen that the universe of $\infty$-groupoids properly contains a universe of sets.
More precisely, there is a subclass of the $\infty$-groupoids of \hott\ which together satisfy the axioms of \etcs.\footnote{In fact, ``\hott'' is not (yet) a single precisely specified theory like \zfc\ and \etcs: as befits a young field, there are many variant theories in use and new ones under development.
  In particular, when I say ``\hott'' I mean to encompass both ``classical'' versions that have the Axiom of Choice and Law of Excluded Middle and also ``intuitionistic'' or ``constructive'' ones that do not.
  In the latter cases, the universe of sets satisfy not \etcs\ (which is classical) but an ``intuitionistic'' version thereof.}
A model of \zfc\ can then be constructed using trees as described in McLarty's chapter, or directly as in~\cite[\S10.5]{hottbook}.
Thus, any mathematics that can be encoded into set theory can also be encoded into \hott.
(Of course, if we intended to encode \emph{all} of mathematics into \hott\ via set theory this way, there would be no benefit to choosing \hott\ as a foundation over set theory.
The point is that \emph{some} parts of mathematics can be also encoded into \hott\ in \emph{other}, perhaps more natural, ways.)

In sum, if we so desire, \emph{we may regard the basic objects of mathematics to be $\infty$-groupoids rather than sets}.
Our discussion in \cref{sec:infty-groupoids} suggests some reasons {why} we might want to do this;
I will mention some further advantages as they arise.
% I will give many different answers to this question in the rest of the chapter, as we explore what \hott\ actually looks like.
% I will return in later sections to the question of , and how such a foundation for mathematics is related to the classical one using set theory.
But it is now time to say something about what \hott\ actually looks like.

\section{Type theory and Logic}
\label{sec:type-theory}

% As the name suggests, homotopy type theory is, in particular, a type theory.
% In this section I will focus mainly on the structure that it shares with all type theories; in the next section I will focus more specifically on its unique features.

The basic objects of \hott\ behave like $\infty$-groupoids; but we generally call them \emph{types} instead, and from now on I will switch to this usage.
This particular word is due to the theory's origins in Martin-L\"{o}f type theory~\cite{martinlof:itt-pred}; but (in addition to being five syllables shorter)
it also fortuitously evokes the terminology ``homotopy type'' from algebraic topology, which is essentially another word for ``$\infty$-groupoid'' (see e.g.~\cite{baez:homotopy-hypothesis}).
 % (although the terminology predated the recognition of this fact by Grothendieck).
%\footnote{Purists will also point out that historically, the collection of homotopy types was actually a set...}

Like sets, the types of \hott\ have \emph{elements}, also called \emph{points}.
We write $x:A$ when $x$ is a point of $A$; the most salient difference between this and \zfc's ``$x\in A$'' is that (like in \etcs) we cannot compare elements of different types: a point is always \emph{a point of some type}, that type being part of its nature.
Whenever we introduce a variable, we must specify its type: whereas in \zfc\ ``for every integer $x$, $x^2\ge 0$'' is shorthand for ``for every thing $x$, if $x$ happens to be an integer then $x^2\ge 0$'', in \hott\ the phrase ``for every integer $x$'' is atomic.
%(I will discuss \hott's logic more carefully in a moment, as it has other distinctive features.)
This arguably matches mathematical practice more closely, although the difference is small.

The basic theory of \hott\ is a collection of \emph{rules} stipulating operations we can perform on types and their points.
For instance, if $A$ and $B$ are types, there is another type called their cartesian product and denoted $A\times B$.
Any such rule for making new types comes with one or more rules for making points of these types: in the case of products, this rule is that given $a:A$ and $b:B$, we have an induced point of $A\times B$ denoted $(a,b)$.
We also have dual rules for extracting information from points of types, e.g.\ from any $x:A\times B$ we can extract $\pi_1(x):A$ and $\pi_2(x):B$.
Of course, $\pi_1(a,b)$ is $a$ and $\pi_2(a,b)$ is $b$.

It is important to understand that these \emph{rules} are not the same sort of thing as the \emph{axioms} of a theory like \zfc\ or \etcs.
Axioms are statements \emph{inside} an ambient superstructure of (usually first-order) logic, whereas the rules of type theory exist at the same level as the deductive system of the logic itself.
In a logic-based theory like \zfc, the ``basic act of mathematics'' is to deduce a conclusion from known facts using one of the rules of logic, with axioms providing the initial ``known facts'' to get started.
By contrast, in a type theory like \hott, the ``basic acts of mathematics'' are specified directly by the rules of the theory, such as the rule for cartesian products which permits us to construct $(x,y)$ once we have $x$ and $y$.
Put differently, choosing the axioms of \zfc\ is like choosing the starting position of a board game whose rules are known in advance, whereas choosing the rules of \hott\ is like choosing the rules of the game itself.
% In particular, unlike \zfc, \hott\ does not depend on first-order logic.
%
% This is a fundamental difference between the above rules for cartesian products and, say, the axiom of pairing in \zfc, which says that for any sets $x$ and $y$, \emph{there exists} a set $z$ such that $\forall w(w\in z \longleftrightarrow (w=x \lor w=y))$.
% Officially, this is a pure \emph{existence} statement, formulated in a background first-order logic.
% By contrast, the rules of \hott\ (like those of any type theory) are \emph{constructions} that directly yield results, rather than merely asserting that objects with certain properties exist; and accordingly, they do not depend on first-order logic.

To understand the effect this distinction has on mathematical practice, we observe that the everyday practice of mathematics can already be separated into two basic activities: constructing (a.k.a.\ defining or specifying) and proving.
For instance, an analyst may first construct a particular function, then prove that it is continuous.
This distinction can be found as far back as Euclid, whose Postulates and Propositions are phrased as things to be \emph{done} (``to draw a circle with any center and radius'') rather than statements of existence, and which are ``demonstrated'' by making a \emph{construction} and then \emph{proving} that it has the desired properties.
Rodin~\cite{rodin:cax} has recently argued that this distinction is closely related to Hilbert's contrast between \emph{genetic} and \emph{axiomatic} methods.\footnote{At least in~\cite{hilbert-bernays:grund-math}; in~\cite{hilbert:zahlbegriff} the same words seem to refer instead to analytic and synthetic theories respectively.}

When encoding mathematics into \zfc, however, the ``construction'' aspect of mathematics gets short shrift, because in fully formal \zfc\
% (as in any theory of first-order logic\footnote{Some presentations of first-order logic permit ``definitional extensions'', but only after it is proven that a unique object with certain properties exists.  For present purposes, this is no different.})
the only thing we \emph{can} do is prove theorems.
Thus, the encoding process must translate constructions into proofs of existence.
By contrast, in \hott\ and other type theories like it, it appears that the pendulum has swung the other way: the \emph{only} thing we can do is perform constructions.
% \footnote{It \emph{is} possible to formulate theories that include both activities on an equal footing; they go by names such as \emph{logic-enriched type theory}~\cite{TODO:aczel}.
%  However, I will argue in a moment that the choice made in \hott, to encode proofs as constructions, is superior.}
How, then, do we encode proofs?

The answer begins with an idea called \emph{propositions as types}: we interpret every \emph{statement} that we might want to prove as a \emph{type}, in such a way that it makes sense to interpret \emph{constructing an element} of that type as \emph{proving the original statement}.
In this way we obtain a form of logic \emph{inside of} type theory, rather than starting with a background logic as is done in set theory.
Thus, as a foundation for mathematics, type theory is ``closer to the bottom'' than set theory: rather than building on the same ``sub-foundations'' (first-order logic), we ``re-excavate'' the sub-foundations and incorporate them into the foundational theory itself.
In the words of Pieter Hofstra, type theory is ``the engine and the fuel all in one.''

One reason this idea is so useful is an observation called the \emph{Curry--Howard correspondence}~\cite{curry:curry-howard,howard:curry-howard,martinlof:itt-pred,wadler:pat}: the logical connectives and quantifiers are \emph{already present} in type theory as {constructions} on types. %\footnote{Surprisingly, the Curry--Howard correspondence (and the propositions-as-types principle) really is originally due to Curry~\cite{curry:curry-howard} and Howard~\cite{howard:curry-howard}, although others such as Martin-L\"{o}f~\cite{martinlof:itt-pred} made important contributions as well; see~\cite{wadler:pat} for an overview.}
For instance, if $A$ and $B$ are types representing propositions $P$ and $Q$ respectively, then $A\times B$ represents the conjunction $P\land Q$.
This is justified because the way we construct an element of $A\times B$ --- by constructing an element of $A$ and an element of $B$ --- corresponds precisely to the way we prove $P\land Q$ --- by proving $P$ and also proving $Q$.
Similarly, the type of functions from $A$ to $B$ (usually denoted $A\to B$) represents the implication $P\to Q$, and so on.

If we interpret logic directly according to this correspondence, we find that just as with the encoding into \zfc, the distinction between construction and proof is destroyed; only this time it is because we have to encode proofs as constructions rather than vice versa.
Whereas in \zfc\ we cannot construct objects, only prove that they exist, under Curry--Howard we cannot prove that something exists without constructing it.

The innovation of \hott\ is to allow both kinds of existence to coexist smoothly.
We follow the overall philosophy of propositions-as-types, but in addition we single out a small but important class of types: those that have at most one point, with no higher equality information.\footnote{The importance of these types has been particularly advocated by Voevodsky, building on precursors such as~\cite{nuprlbook,ab:bracket-types}.}
I will call these types \emph{truth values}, since we think of them as representing ``false'' (if empty) or ``true'' (if inhabited); they are also often called \emph{propositions} or \emph{mere propositions}.
Moreover, we add a rule that for any type $A$ there is a \emph{truncation} $\merely{A}$ (also called the \emph{bracket} or \emph{squash}), such that $\merely A$ is a truth value, and such that given any $a:A$ we have $\tr a : \merely A$.
(Since $\merely A$ is a truth value, $\tr a$ doesn't depend on the value of $a$, only that we have it.)

Now we can distinguish between existence proofs and constructions by whether the type of the result is truncated or not.
When we construct an element of a type $A$ that is not a truth value, we are defining some specific object; but if we instead construct an element of $\merely A$, we are ``proving'' that some element of $A$ exists without specifying it.\footnote{The possibility of these two interpretations of existence was actually already noticed by Howard~\cite[\S12]{howard:curry-howard}.}
From this point of view, which is shared by many members of the \hott\ community, it is misleading to think of propositions-as-types as ``encoding first-order logic in type theory''.
While this description can serve as a first approximation, it leads one to ask and argue about questions like ``should the statement $\exists x:A$ be encoded by the type $A$ or the type $\merely A$?''
We regard this question as invalid, because it implicitly assumes that mathematics has already been encoded into first-order logic, with constructions and pure-existence proofs collapsed into the quantifier $\exists$.
We reject this assumption: the proper approach is to encode \emph{mathematics} directly into \hott, representing a construction of an element of $A$ by the type $A$ itself, and a pure-existence statement by its truncation $\merely A$.

It is true that due to the ascendancy of \zfc\ and first-order logic in general, most modern mathematicians ``think in first-order logic'' and are not used to distinguishing constructions from existence proofs.
However, it remains true that some kinds of theorem, such as ``$A$ is isomorphic to $B$'', are almost always ``proven'' by giving a construction; and a careful analysis reveals that such ``proofs'' have to convey more information than mere existence, because frequently one needs to know later on exactly \emph{what} isomorphism was constructed.
%(This is the cause of many uses of weasel words such as ``canonical''.)
This is one of the ways in which \hott\ represents the actual practice of mathematics more faithfully than other contenders.
With a little bit of practice, and careful use of language, we can learn to consciously use this feature when doing mathematics based on \hott.
 %\footnote{The novelty of this approach does mean that we need careful conventions of language to use it correctly, and the \hott\ community has not yet reached a clear consensus on what these conventions should be.
% The first edition of~\cite{hottbook} used ``there exists'' for constructions and ``there merely exists'' for proofs, but recently momentum seems to be gathering around Voevodsky's proposal of ``there exists'' for proofs and something like ``there is a given'' for constructions.}

By the way, while the distinction between construction and proof is sometimes identified with the opposition between constructive/intuitionistic and classical logic (as is suggested by the shared root ``construct''), the relationship between the two is actually limited.
On one hand, while it is true that the ``natural'' logic obtained by Curry--Howard turns out to be intuitionistic, one can add additional axioms %\footnote{Under the type-theoretic view of logic, an ``axiom'' is essentially a rule saying that some specified type has an element, even if the standard rules would not give us a way to construct one.}
that are not ``constructive'' but can nevertheless be used in ``constructions''.
%(which to avoid confusion it might be better to call ``specifying'' them).
Indeed, the exceedingly nonconstructive Axiom of Choice asserts exactly that objects which merely exist can nevertheless be assumed to be specified,
i.e. ``constructed'' in a formal sense.
In particular, axioms of classical logic can consistently be included in \hott.

On the other hand, intuitionistic first-order logic includes ``pure unspecified existence'' just like classical logic does, and constructive/intuitionistic set theory~\cite{beeson:fcm,aczel:cst} collapses constructions into proofs just like \zfc\ does.
It is true that constructive mathematicians in the tradition of Martin-L\"{o}f~\cite{martinlof:itt-pred} do adhere intentionally to the original Curry--Howard interpretation, % (we cannot prove that something exists without defining it),
regarding it as part of their constructivism; but they must also separately refrain from using any nonconstructive principles.
That is, a constructive \emph{philosophy} may lead one to prefer ``constructions'' to proofs, but this is a separate thing from the (intuitionistic) \emph{logic} that it also leads one to prefer.
Moreover, Escard\'o has recently argued that Brouwer himself must have intended some notion of unspecified existence, since his famous theorem that all functions $\mathbb{N}^{\mathbb{N}}\to \mathbb{N}$ are continuous is actually \emph{inconsistent} under unmodified Curry--Howard~\cite{escardo-xu:brouwer-ch}.

A last aspect of type theory that deserves mention is its computational character: its rules can also be read as defining a programming language that can actually be executed by a computer.
This makes it especially convenient for computer formalization of mathematical proofs, as well as for mathematical verification of computer software.
Thus, \hott\ is also better-adapted to these purposes than set theory is,\footnote{Although finding the best way to extend the computational aspects of type theory to the specific features of \hott\ is an active research area.} and indeed computer formalization has been very significant in the origins and development of \hott.
But this would fill a whole chapter by itself, so reluctantly I will say no more about it here.

\section{Identifications and equivalences}
\label{sec:ident-equiv}
\label{sec:cardinality}

So far, I have not really said anything that is unique to \hott.
The description of types, rules, and elements in \cref{sec:type-theory} applies to any type theory, including Martin-L\"{o}f's original one.
The approach to logic using truncations is more novel, but it still does not depend on regarding types as $\infty$-groupoids.
%; it could also be pursued in a more traditional type theory whose types are thought of as mere sets.
However, this kind of logic is particularly appropriate in \hott, for several reasons.

The first is that, like our considerations in \cref{sec:infty-groupoids}, it drives us inexorably from sets to % groupoids and
$\infty$-groupoids.
Namely, if statements are interpreted by types, then in particular for any $x:A$ and $y:A$, the statement ``$x=y$'' must be \emph{a type}, whose points we refer to as \emph{identifications} of $x$ with $y$.
If $A$ is a set, then this type is a mere truth value, but in general there is no reason for it to be so.
%; there could be many distinct indentifications of $x$ with $y$.

Somewhat magically, it turns out that the ``most natural'' rule governing the type $x=y$, as first given by Martin-L\"of~\cite{martinlof:itt-pred}, does \emph{not} imply that it is always a truth value, but \emph{does} imply that it automatically inherits the structure of an $\infty$-groupoid~\cite{pll:wkom-type,bg:type-wkom}.
This rule is related to Leibniz's ``indiscernibility of identicals'', but its form is rather that of Lawvere~\cite{lawvere:comprehension}, who characterized equality using an adjunction between unary predicates and binary relations.
Martin-L\"of's version says that if we have a type family $C(x,y,p)$ depending on $x$, $y$, and an identification $p:x=y$, and if we want to construct an element of $C(x,y,p)$ for every $x$, $y$, and $p$, then it suffices to construct elements of $C(x,x,\mathsf{refl}_x)$ for every $x$.
(Here $\mathsf{refl}_x$ denotes a canonically specified element of $x=x$, called the \emph{reflexivity witness} or the \emph{identity identification}.
The standard proofs of transitivity and symmetry of equality from indiscernibility of identicals become in \hott\ constructions of the first level of $\infty$-groupoid structure.)

% (The classical rule of indiscernibility of identicals, as described in Baez's chapter, is the restriction of this rule to the special case when we are proving rather than constructing, and when the statement to be proven is independent of the witness $p$.
% Baez mentions how this this rule implies transitivity and symmetry of equality; in \hott\ those implications becomes a construction of the first level of $\infty$-groupoid structure.)

In this way, $\infty$-groupoids become much simpler in \hott\ than they are in set theory.
We saw in \cref{sec:infty-groupoids} that a mathematician trying to \emph{define} $\infty$-groupoids in set theory is led to a rather complicated structure. %; indeed, from that point of view, one would be quite justified in questioning whether such beasts could really be considered ``foundational'' (even though, as we have seen, they do satisfy the formal notion of ``foundation'' that I described in \cref{sec:synthetic}).
However, \hott\ reveals that \emph{synthetically}, an $\infty$-groupoid is really quite a simple thing: we might say that we obtain a synthetic theory of $\infty$-groupoids by (0) starting with type theory, (1) taking seriously the idea that a statement of equality $x=y$ should be a type, (2) writing down the most natural rule governing such a type, and then (3) simply \emph{declining to assert} that all such types are mere truth values.\footnote{Many people contributed to this view of Martin-L\"{o}f's equality types, but Hofmann--Streicher~\cite{hs:gpd-typethy} and Awodey--Warren~\cite{aw:htpy-idtype} were significant milestones.}

To be precise, however, this is not quite correct; better would be to say that Martin-L\"{o}f's type theory, unlike set theory, is sufficiently general to \emph{permit} its types to be treated as $\infty$-groupoids, and in \hott\ we choose to do so.
(This is analogous to how intuitionistic logic, unlike classical logic, is sufficiently general to permit the assumption of topological or computational structure.)
%, leading to ``synthetic topology''~\cite{bl:met-syntop} or ``synthetic computability theory''~\cite{bauer:syncomp,richman:churchs-thesis}.
Thus, in order to obtain a true synthetic theory of $\infty$-groupoids, we need to add some rules that are \emph{specific} to them, which in particular will ensure that it is definitely \emph{not} the case that all equality types are truth values.

The principal such rule in use is Voevodsky's \emph{univalence axiom}~\cite{klv:ssetmodel}.
This is formulated with reference to a \emph{universe type} $\U$, whose points are other types.
(For consistency, $\U$ cannot be a point of itself; thus one generally assumes an infinite hierarchy of such universes.)
Univalence says that for types $A:\U$ and $B:\U$, the type $A=B$ consists of \emph{equivalences} between $A$ and $B$, the latter being a standard definition imported from higher category theory\footnote{Although it requires some cleverness to formulate it correctly in type theory; this was first done by Voevodsky.} that generalizes bijections between sets.
In particular, if $A$ has any nontrivial automorphisms, then $A=A$ is not a mere truth value. %, then $A=A$ has more than one distinct element.

Univalence is the central topic of Awodey's chapter; he concludes that it codifies exactly the principle of structuralism, ``isomorphic objects are identical''.
Indeed, with univalence we no longer need any Fregean abstraction to define ``structure''; we can simply consider types themselves (or, more generally, types equipped with extra data) to \emph{be} structures.
Fregean abstraction is for forgetting irrelevant facts not preserved by isomorphism, like whether $0\in 1$; but in \hott\ there are no such facts, since isomorphic types are actually already \emph{the same}.
Thus, if we wish, we may consider \hott\ to be a \emph{synthetic theory of structures}.\footnote{This is not in conflict with also calling it a synthetic theory of $\infty$-groupoids; the two phrases simply emphasize different aspects of \hott.
We could emphasize both aspects at once by calling it a ``synthetic theory of $\infty$-groupoidal structures''.}

More concretely, univalence ensures that any construction or proof can be transported across an isomorphism (or equivalence): anything we prove about a type is automatically also true about any equivalent type.
Here again \hott\ captures precisely an aspect of mathematical practice that is often glossed over by set theory.
Univalence also implies that two ``truth values'', as defined in \cref{sec:type-theory},
are equal as soon as they are logically equivalent; thus they really do carry no more information than a truth value.

A second way that the logic of \cref{sec:type-theory} is particularly appropriate is that \hott\ clarifies the distinction between types and truth values, by placing it on the first rung of an infinite ladder.
In fact, for any integer $n\ge -2$ there is a class of types called \emph{$n$-types},\footnote{This notion is well-known in homotopy theory under the name \emph{homotopy $n$-type} and in higher category theory under the name \emph{$n$-groupoid}.  Its definition in type theory is due to Voevodsky, who calls them ``types of h-level $n+2$''.} % instead of $n$-types, thereby starting the numbering at $0$ instead of $-2$.  This is arguably more logical, but contrary to decades of convention.}
  such that the singleton is the only $(-2)$-type, the truth values are the $(-1)$-types, and the sets are the $0$-types.
Informally, an $n$-type contains no higher equality information above level $n$: two elements of a $0$-type (i.e.\ a set) can be equal in at most one way, two \emph{equalities} in a $1$-type can be equal in at most one way, and so on.
Formally, $A$ is an $n$-type if for all $x:A$ and $y:A$, the type $x=y$ is an $(n-1)$-type (with the induction bottoming out at $n=-2$).

In addition, for any $n$ we have an \emph{$n$-truncation} operation: $\trunc{n}{A}$ is an $n$-type obtained from $A$ by discarding all distinctions between equalities above level $n$. %, obtaining an $n$-type in a universal way.
In particular, $\trunc{-1}{A}$ discards all distinctions between \emph{points} of $A$, remembering only whether $A$ is inhabited; thus it is the truth-value truncation $\merely{A}$ from \cref{sec:type-theory}.
The next most important case is the $0$-truncation $\trunc{0}{A}$, which makes $A$ into a set by discarding distinctions between equalities between its points, remembering only the truth value of whether or not they are equal.

At this point we can deal with one of the examples of a groupoid from \cref{sec:infty-groupoids}: sets and cardinalities.
In \hott\ the \emph{type of sets} is naturally defined as a subtype of the universe $\U$ which contains only the sets ($0$-types).
By univalence, then, for sets $A$ and $B$, the type $A=B$ is the type of bijections between them.
Thus two sets are \emph{automatically} identical exactly when they are bijective, so it may appear that there is no need to specify the equalities separately from the points in this case.

However, since the type of bijections between sets $A$ and $B$ is itself a set and not (generally) a truth value, the type of sets is a $1$-type and not a set. %\footnote{More generally, the type of $n$-types is an $(n+1)$-type and not an $n$-type, and the full universe $\U$ is not an $n$-type for any finite $n$.}
This is an important difference with \zfc, in which the collection of sets (belonging to some universe) \emph{is} itself a set --- but it matters little in mathematical practice, which is mostly structural.
Indeed, mathematicians familiar with category theory tend to be drawn to this idea: it seems perverse to distinguish between isomorphic sets as \zfc\ does.\footnote{This should not be confused with distinguishing between \emph{subsets} of some fixed set that may be abstractly isomorphic as sets, such as $\mathbb{N}\subseteq\mathbb{R}$ and $\mathbb{Q}\subseteq\mathbb{R}$, which is common and essential to mathematics.
The point is rather that of Benacerraf~\cite{benacerraf:wncnb}: there is no reason to distinguish between, say, $\{\emptyset, \{\emptyset\}, \{\{\emptyset\}\}, \dots \}$ and $\{\emptyset, \{\emptyset\}, \{\emptyset,\{\emptyset\}\}, \dots \}$ as definitions of ``the natural numbers''.}
% , based on a complicated hierarchical membership predicate.

On the other hand, mathematicians \emph{are} accustomed to consider the collection of \emph{cardinalities} to form a set (modulo size considerations).
Thus, in \hott\ it is sensible to define the set of {cardinalities} to be the $0$-truncation of the type of sets.
That is, a cardinality is presented by a set, and bijective sets present equal cardinalities; but unlike sets, two cardinalities can be equal in at most one way.
One nice consequence is that the subset of \emph{finite} cardinalities is then equal to the natural numbers.

% We have not yet discussed how Bishop's general set-definition philosophy is to be rendered in \hott, but the $0$-truncation can certainly be seen as a particular instance of it: given $A$, the set $\trunc 0 A$ is defined by (1) to construct a point of $\trunc 0 A$, we give a point of $A$, and (2) to show that $x,y:A$ yield equal elements of $\trunc 0 A$, we show that they are equal in $A$.
% Thus, our definition of the set of cardinalities is also included.
% In contrast to \cref{sec:infty-groupoids}, however, we did not have to explicitly insert the equivalence relation of equinumerosity; we automatically obtained isomorphism-invariance from the univalence axiom.

The $0$-truncation has many other uses; for instance, it allows us to import the definition of \emph{homotopy groups} from algebraic topology.
Given a type $X$ and a point $x:X$, we first define the \emph{loop space} $\Omega(X,x)$ to be the type $x=x$, and the \emph{$n$-fold loop space} by induction as $\Omega^{n+1}(X,x) = \Omega^n(\Omega(X,x),\mathsf{refl}_x)$.
The \emph{$n^{\mathrm{th}}$ homotopy group} of $X$ based at $x$ is then $\pi_n(X,x) = \trunc 0 {\Omega^n (X,x)}$.
% This is, by definition, a set, which has the structure of a group induced by the transitivity of identifications.
If $X$ is an $n$-type, then $\pi_k(X,x)$ is trivial whenever $k>n$; in general it can be said to measure the nontriviality of the identification structure of $X$ at level $n$.
For instance, if $X$ is a set, then $\pi_k(X,x)$ is trivial for any $k\ge 1$; whereas if $X=\U$ and $x$ is a set $A$, then $\pi_1(\U,A)$ is the automorphism group of $A$ while $\pi_k(\U,A)$ is trivial for $k>1$.

\section{Higher inductive types}
\label{sec:high-induct-types}

I mentioned in \cref{sec:type-theory} that \hott\ consists of rules describing operations we can perform on types and their points.
In fact, all but a couple of these rules belong to one uniformly specified class, known as \emph{higher inductive types} (\hits), which can be considered a generalization of Bishop's rule for set-construction that takes higher identifications into account.

Higher inductive types include, in particular, \emph{ordinary} inductive types, which have been well-known in type theory for a long time (several examples appear already in~\cite{martinlof:itt-pred}).
The simplest sorts of these are \emph{nonrecursive}, in which case the rule says that to define a type $X$, we specify zero or more ways to construct elements of $X$.
This amounts to stipulating some finite list of functions with codomain $X$ and some specified domain, called the \emph{constructors} of $X$.
For instance, given types $A$ and $B$, their \emph{disjoint union} $A+B$ is specified by saying that there are two ways to construct elements of $A+B$, namely by injecting an element of $A$ or an element of $B$; thus we have two constructors $\mathsf{inl}:A\to A+B$ and $\mathsf{inr}:B\to A+B$.

As recognized in~\cite[\S1.1]{martinlof:itt-pred}, this is similar to Bishop's rule; the main difference is that we omit the specifying of equalities.
How then are we to know when two points of such a type are equal?
The answer is that an inductive type should be regarded as \emph{freely generated} by its constructors, in the sense that we do not ``put in'' anything --- whether a point or an identification --- that is not \emph{forced} to be there by the constructors.
For instance, every point of $A+B$ is either of the form $\mathsf{inl}(a)$ or $\mathsf{inr}(b)$, since the constructors do not force any other points to exist.
Moreover, no point of the form $\mathsf{inl}(a)$ is equal to one of the form $\mathsf{inr}(b)$, since the constructors do not force any such identifications to exist.
However, if we have $a:A$ and $a':A$ with $a=a'$, then there \emph{is} an induced identification $\mathsf{inl}(a)=\mathsf{inl}(a')$, since all functions (including $\mathsf{inl}$) must respect equality.
% Along these lines, one can show that $A+B$ is a set precisely when both $A$ and $B$ are.

More generally, ordinary inductive types can be \emph{recursive}, meaning that some of the constructors of $X$ can take as input one or more {previously} constructed elements of $X$.
% These are called \emph{recursive} constructors.
For example, the natural numbers $\mathbb{N}$ have one nonrecursive constructor $0:\mathbb{N}$ and one recursive one $s:\mathbb{N}\to \mathbb{N}$.
The elements and equalities in such a type are all those that can be obtained by applying the constructors, over and over again if necessary.

\emph{Higher} inductive types are a generalization of ordinary ones, which were invented by the author and others.\footnote{Specifically, Lumsdaine, Bauer, and Warren, with further notable contributions by Brunerie and Licata.
  The basic theory of \hits\ is still under development by many people; currently the best general reference is~\cite[Ch.~6]{hottbook}.}
The simplest case is a nonrecursive \emph{level-1} \hit, where in addition to specifying ways to construct elements of $X$, we can specify ways to construct identifications between such elements.
Thus, in addition to constructor functions as before (which we now call \emph{point-constructors}), we also have \emph{identification-constructors}.

This is almost the same as Bishop's rule for set-construction, with two differences.
Firstly, a \hit\ need not be a set. % --- though we can always $0$-truncate it.
Secondly, the identification-constructors need not form an equivalence relation; e.g. we may specify $x=y$ and $y=z$ but not $x=z$.
However, since all types \emph{are} $\infty$-groupoids,
% the ``free generation'' rule that we put in everything we are forced to (and nothing else) demands that 
in such a case it will nevertheless be \emph{true} that $x=z$.
More precisely, if we have constructors yielding identifications $p:x=y$ and $q:y=z$, then there will be an induced identification $p\ct q : x=z$, which is forced to exist even though we didn't ``put it in by hand''.

Suppose now that we \emph{are} in Bishop's situation, i.e.\ we have a type $A$ and an equivalence relation $\sim$ on it.
We can define a \hit\ $X$, with one point-constructor $q:A\to X$, and one identification-constructor saying that whenever $a\sim a'$ we have $q(a) = q(a')$.
Then $X$ will be close to the quotient of $\sim$, except that it will not generally be a set even if $A$ is.
For instance, since $a\sim a$ for any $a:A$, our identification-constructor yields an identification $q(a)=q(a)$; but nothing we have put into $X$ forces this identification to be the same as $\mathsf{refl}_{q(a)}$, and so (by the free generation principle) it is not.
Thus, to obtain the usual quotient of $\sim$, we have to $0$-truncate $X$; in \hott\ we may call this the \emph{set-quotient}.
For instance, the set of real numbers could be defined as the set-quotient of the equivalence relation on infinite decimal expansions from \cref{sec:infty-groupoids}.
% The real numbers $q(0.5)$ and $q(4.\overline{9})$ will then really be \emph{the same}.
In this way we essentially recover Bishop's set-formation rule.\footnote{There is one subtle difference: Bishop actually allows us to distinguish between $0.5$ and $0.4\overline{9}$ as long as we speak of an ``operation'' rather than a ``function''.
  In \hott\ such an ``operation'' is just a function defined on decimal expansions, not anything acting on ``real numbers''.}
  % In \hott, the real numbers $q(0.5)$ and $q(0.4\overline{9})$ are honestly equal, hence indistinguishable by any construction we can perform.
  % We do, of course, work with real numbers by representing them as decimal expansions; but the rule governing the set-quotient says that for a construction to induce an operation on \emph{the set of real numbers}, it must respect such equalities.
  % Bishop instead considered a quotient to still ``remember'' the type $A$, allowing constructions on ``real numbers'' to distinguish between $0.5$ and $0.4\overline{9}$ (though he denigrated such constructions by calling them ``operations'' rather than ``functions'').
  % We can perform such constructions in \hott\ as well, by working with decimal expansions directly, but we regard it as misleading to claim that they act on ``real numbers''.
  % Indeed, essentially this point apparently misled Bishop himself once into believing that the axiom of choice should be automatic (see~\cite[Chapter 1 Notes]{bb:constr-analysis}).}

Higher inductive types can also be recursive: both kinds of constructor can take previously constructed elements of $X$ as inputs.
This is very useful --- e.g.\ it yields free algebraic structures, homotopical localizations, and even the $n$-truncation --- but also somewhat technical, so I will say no more about it.

The reader may naturally wonder \emph{why} we don't ask the identification-constructors to form an equivalence relation. %, rather than supposing all the $\infty$-groupoid structure to be freely generated.
One reason is that for \hits\ that are not sets, the analogue of an equivalence relation would be an ``$\infty$-groupoid'' in the exceedingly complicated sense referenced at the beginning of \cref{sec:synthetic}.
Forcing ourselves to use such structures would vitiate the already-noted advantages of a \emph{synthetic} theory of $\infty$-groupoids.

As a concrete example of the usefulness of {not} requiring equivalence relations \emph{a priori}, if we have two functions $f,g:A\rightrightarrows B$ between sets, we can construct their \emph{set-coequalizer} as the $0$-truncation of the \hit\ with one point-constructor $q:B\to X$ and one identification-constructor saying that for any $a:A$ we have $q(f(a))=q(g(a))$.
In set theory, we would have to first construct the equivalence relation on $B$ freely generated by the relations $f(a)\sim g(a)$ and then take its quotient; \hits\ automate that process for us.
Moreover, if we omit the assumption that $A$ and $B$ are sets and also omit the $0$-truncation, we obtain a \emph{homotopy coequalizer}, which would be \emph{much} harder to construct otherwise.

Another reason for considering freely generated $\infty$-groupoids is that many very interesting $\infty$-groupoids \emph{are} freely generated, and in most cases a fully explicit description of them \emph{is not known} and is not expected to be knowable.
Thus, \hits\ are the \emph{only} way we can represent them in \hott.

A simple example of a freely generated $\infty$-groupoid is the \emph{circle}\footnote{This is a ``homotopical'' circle, not a ``topological'' circle such as $\{(x,y)\in \mathbb{R}\times \mathbb{R} \mid x^2+y^2=1\}$.
  The latter can also be defined in \hott, of course, but it will be a set, whereas the \hit\ $\mathbb{S}^1$ is not.
  The homotopical circle is so-called because it is the \emph{shape} (a.k.a. ``fundamental $\infty$-groupoid'') of the topological circle, with continuous paths in the latter becoming identifications in the former; and historically $\infty$-groupoids were originally studied as shapes of topological spaces.
  In \hott\ the shape ought to be constructible as a \hit, but no one has yet managed to do it coherently at all levels.
  Unlike classically, not every type in \hott\ can be the shape of some space, but we can hope that the \hit\ $\mathbb{S}^1$ is still the shape of the topological circle.

  There is an arguably better approach to such questions called ``axiomatic cohesion''~\cite{ss:qgftchtt,shulman:bfp-realcohesion}, in which the types of \hott\ are enhanced to carry intrinsic topological structure in addition to their higher identifications.
  Unfortunately, space does not permit me to discuss this here, but a brief introduction can be found in Corfield's chapter.\label{fn:cohesion}}
$\mathbb{S}^1$, which as a \hit\ has one point-constructor $\mathsf{b}:\mathbb{S}^1$ and one identification-constructor $\ell:\mathsf{b}=\mathsf{b}$.
Since nothing forces $\ell$ to be equal to $\mathsf{refl}_{\mathsf{b}}$, it is not ---
nor is $\ell\ct\ell$, or $\ell \ct \ell\ct \ell$, and so on.
In fact, $\Omega(\mathbb{S}^1,\mathsf{b})$ is isomorphic to the integers $\mathbb{Z}$.\footnote{This is well-known in homotopy theory; its first proof in \hott\ by the author~\cite{ls:pi1s1} was an early milestone in combining \hits\ with univalence.}
Since $\mathbb{Z}$ is a set, this implies $\pi_1(\mathbb{S}^1,\mathsf{b})=\mathbb{Z}$ while $\pi_k(\mathbb{S}^1,\mathsf{b})=0$ for all $k>1$, so in this case we do have a fully explicit description.
However, there are similar types for which no such characterization is known, particularly when we move on to \emph{level-$n$} \hits\ having constructors of ``higher identifications''.
For instance, the \emph{2-sphere} $\mathbb{S}^2$ has one point-constructor $\mathsf{b}:\mathbb{S}^2$ and one level-2 identification-constructor $\mathsf{refl}_{\mathsf{b}} = \mathsf{refl}_{\mathsf{b}}$; the \emph{3-sphere} has $\mathsf{b}:S^3$ with a level-3 $\mathsf{refl}_{\mathsf{refl}_{\mathsf{b}}}= \mathsf{refl}_{\mathsf{refl}_{\mathsf{b}}}$; and so on.
Analogously to $\mathbb{S}^1$ we have $\pi_n(\mathbb{S}^n)=\mathbb{Z}$,\footnote{Also a standard result in homotopy theory; see~\cite[Ch.~8]{hottbook} and~\cite{lb:pinsn} for proofs in \hott.} but also for example $\pi_3(\mathbb{S}^2)=\mathbb{Z}$, despite the fact that $\mathbb{S}^2$ has no \emph{constructors} of level 3.
In general, $\pi_k(\mathbb{S}^n)$ is usually nontrivial when $k\ge n$, but most of its values are not known.
Computing them, for classically defined $\infty$-groupoids, is a major research area which is not expected to ever be ``complete''.

What does this mean to a philosopher?
For one thing, it shows how a simple foundational system can give rise very quickly to deep mathematics.
The rules governing \hits\ are arguably unavoidable, once we have the idea of defining types in such a way; % even more so than the axioms of \zfc,
% \footnote{Since there isn't space to give those rules explicitly here, you'll have to either take my word for this or go read about them elsewhere, e.g. in~\cite{hottbook}.}
while the spheres $\mathbb{S}^n$ result from quite simple applications of those rules.
Moreover, we have seen % (in \cref{sec:infty-groupoids} and Baez's chapter)
that even the basic notion of $\infty$-groupoid arises inescapably from thinking about {equality}.
Thus, there are {numerical} invariants like $\pi_3(\mathbb{S}^2)$ % (the cardinalities of the $\pi_k(\mathbb{S}^n)$, most of which are finite)
quite close to the foundations of logic.

% \section{On models}
% \label{sec:models}

% The discussion of homotopy groups of spheres at the end of the last section brings into focus a question which I have not carefully addressed yet: how is the synthetic theory of $\infty$-groupoids in \hott\ related to the analytic theory of $\infty$-groupoids as studied in classical mathematics?

% % 

% (So far, in \hott\ we have computed a few of these groups and found them to be the same as the classical ones; a major open question is how many of their values are determined.
% In some limited ways, this is analogous to the Continuum Problem of set theory: the cardinality of $\mathbb{R}$ can also be considered a ``numerical invariant'' (though an infinite one), which G\"{o}del eventually showed was \emph{not} determined by the axioms of \zfc.
% Here, of course, we have a somewhat different situation in that many of the homotopy groups of spheres have known values from classical mathematics, and the question is just whether \hott\ is sufficient to determine \emph{that} value.)

\section{General Covariance}
\label{sec:erlangen}

At long last, we return to the third example from \cref{sec:infty-groupoids}: spacetime manifolds.
For simplicity, I will consider only \emph{Minkowski} spacetimes, corresponding to special rather than general relativity; similar ideas can be applied to other kinds of gauge invariance/covariance as well.

A modern mathematician defines a Minkowski spacetime to be a 4-dimensional real affine space with a Lorentzian inner product.
We can repeat this definition in \hott, yielding a type $\mink$ whose points are Minkowski spacetimes.
Now we can ask what the \emph{identifications} are in $\mink$.
This is a special case of a more general question: what are the identifications in a \emph{type of structured sets}?
Recall that univalence ensures that identifications in the type of \emph{all} sets are bijections; this turns out to imply that an identification of structured sets is a bijection which ``preserves all the structure'', i.e.\ an \emph{isomorphism} in the appropriate category (see e.g.~\cite[\S9.8]{hottbook}).
Thus, % an identification in the type of groups is a group isomorphism, an identification in the type of manifolds is a diffeomorphism --- and
an identification in $\mink$ is an isometry, as we would hope.
In particular, anything we can say in \hott\ about Minkowski spacetimes is automatically covariant under isometry.

Note that since isometries form a set, $\mink$ is a 1-type.
We could, if we wished, $0$-truncate it to obtain a set, as we did with the type of sets in \cref{sec:cardinality} to obtain the set of cardinalities.
However, the hole argument tells us that this would be \emph{wrong}, at least for the purpose of modeling reality: we really do need to remember the nontrivial identifications in $\mink$.

So far, so good.
However, there is another side to the story, which I alluded to briefly in \cref{sec:infty-groupoids}: why did the hole argument seem paradoxical for so long?
This can be attributed at least partly to a radically different viewpoint on manifolds, as described by Norton~\cite{norton:covariance}:
\begin{quotation}
%  There is a presumption in much modern interpretation of Einstein\dots that much of what he says cannot be taken at face value.  (Why does Einstein make such a fuss about introducing arbitrary spacetime coordinates?  We have always been able to label spacetime events any way we please!)\dots  [My] proposal\dots is that
  \dots our modern difficulty in reading Einstein literally actually stems from a change\dots in the mathematical tools used\dots.
  In recent work\dots we begin with a very refined mathematical entity, an abstract differentiable manifold\dots.  We then judiciously add further geometric objects only as the physical content of the theory warrants\dots.
  In the 1910s, mathematical practices in physics were different\dots.  one used number manifolds --- $\mathbb{R}^n$ or $\mathbb{C}^n$ for example.  Thus Minkowski's `world'\dots was literally $\mathbb{R}^4$, that is it was the set of all quadruples of real numbers.

  Now anyone seeking to build a spacetime theory with these mathematical tools of the 1910s faces very different problems from the ones we see now.  Modern differentiable manifolds have too little structure and we must add to them.  Number manifolds have far too much structure\dots the origin $\langle 0,0,0,0\rangle$ is quite different from any other point, for example\dots. The problem was not how to add structure to the manifolds, but how to deny physical significance to existing parts of the number manifolds.  How do we rule out the idea that $\langle 0,0,0,0\rangle$ represents the preferred center of the universe\dots?
\end{quotation}
In brief, \emph{mathematical structuralism} had not yet been invented.
Our explanation of the hole argument relied on comfort with the structural idea of an isometry between abstract manifolds. % we were then able to realize that the identity map is not special, and in particular cannot be used as an identification if it does not preserve the metrics.
But if one views spacetime as the \emph{specific} manifold $\mathbb{R}^4$, this sort of argument is unavailable; thus the confusion surrounding the hole argument becomes more understandable.

While structuralism is the modern method of choice to deal with this conundrum, it is not the only possible solution;
historically, Klein's \textit{Erlangen} program was used for the same purpose.
Here is Norton again:
\begin{quote}
  Felix Klein's \textit{Erlangen} program provided precisely the tool that was needed.  One assigns a characteristic group to the theory\dots.  Only those aspects of the number manifold that remain invariant under this group are allowed physical significance\dots.  As one increases the size of the group, one strips more and more physical significance out of the number manifold.
\end{quote}
This suggests a different definition of $\mink$:
%According to the $\infty$-groupoidified reading of Bishop incarnated by \hits, when defining a type we are allowed to specify not only its points, but also the identifications between them.
%If we define $\mink$ in a structuralist way as above, with its points being abstract Lorentzian affine spaces, then univalence ensures that all the desired identifications between these points are already present.
% However, with the \textit{Erlangen} viewpoint in mind,
we could begin with the singleton type $\{\mathbb{R}^4\}$ and \emph{add identification-constructors} making up the desired symmetry group (in this case, the Poincar\'{e} group\footnote{The Poincar\'{e} group is usually considered not as a discrete group but as a \emph{Lie} group, with its own manifold structure.
  This can be incorporated as well using ``axiomatic cohesion'', mentioned briefly in \cref{fn:cohesion} on page~\pageref{fn:cohesion}.}). % see~\cite{lf:emspaces} for the precise construction.
%, consisting of the automorphisms of $\mathbb{R}^4$ that are generated by translations, 3D rotations, and Lorentz boosts).
% This is a special case of the construction of Eilenberg--MacLane spaces as in~\cite{lf:emspaces}.
In other words, % rather than defining a Minkowski spacetime to be an arbitrary set equipped with various kinds of structure,
we say that there is \emph{one} Minkowski spacetime, namely $\mathbb{R}^4$, and that it can be identified with itself in many ways, such as translations, 3D rotations, and Lorentz boosts.
These extra added identifications force everything we say about ``Minkowski spacetimes'' to be invariant under their action.
For example, while in $\mathbb{R}^4$ we can distinguish the point $\langle 0,0,0,0\rangle$, in a {Minkowski spacetime} we cannot, because this point is not invariant under translations.
However, we can say that a Minkowski spacetime comes with a Lorentzian distance function, since this structure on $\mathbb{R}^4$ \emph{is} preserved by the Poincar\'{e} group.
This is precisely the point of the \textit{Erlangen} program, %: only the ``geometric'' notions are invariant under the chosen group action.
which \hott\ codifies into the foundations of mathematics by constructing a type that ``remembers exactly those aspects of $\mathbb{R}^4$ preserved by the group action.''

Finally, we can show in \hott\ that these two definitions of Minkowski spacetime agree.
Roughly, this is because two abstract Minkowski spacetimes can always be identified \emph{somehow}, while their automorphisms can be identified with the Poincar\'{e} group; thus the points and the identifications can be matched up consistently.
Thus, \hott\ could be said to unify the \textit{Erlangen} and structuralist approaches to geometry.

One might argue that these approaches were unified long ago, by the development of category theory.
Indeed, as detailed in~\cite{marquis:gpov}, category theory can be seen as a generalization of the \textit{Erlangen} program, where rather than simply having a group act by automorphisms of a single object, we consider isomorphisms, or more generally morphisms, between different objects, and permit as meaningful only those properties that vary appropriately under such transformations (i.e.\ those that are covariant --- or, perhaps, \emph{contravariant}, the dual sort of variation that can be distinguished only once we allow noninvertible morphisms).
And category theory is, of course, the language of choice for the modern structuralist.

However, when category theory is built on top of a foundational set theory, one has to take the additional step of \emph{defining} the notion of isomorphism as the appropriate ``criterion of sameness'' and (in principle) \emph{proving} that all properties of interest are invariant under isomorphism.
As Marquis says, in the \textit{Erlangen} program:
\begin{quote}
  \dots what is usually taken as a \emph{logical} notion, namely equality of objects, is captured in geometry by motions, or transformations of the given group. \cite[p19; emphasis added]{marquis:gpov}
\end{quote}
Moreover, when generalized to higher groupoids and higher categories, this leads to the highly complicated \emph{defined} notion of $\infty$-groupoid mentioned in \cref{sec:synthetic}.
But with univalence and \hits, \hott\ places the notion of equality back where it belongs --- in logic, or more generally the foundations of mathematics --- while maintaining the insights of the \textit{Erlangen} program and category theory.

% Note that the idea underlying Bishop's sets as well as \hits{} --- that when specifying a set or type we must specify its identifications as well as its elements --- is already implicitly present in the \textit{Erlangen} program:
% the same space $\mathbb{R}^4$ becomes a different geometry based on whether we identify it with itself by Lorentz transformations or Euclidean transformations.
% Moreover, as pointed out in Corfield's chapter, Cassirer~\cite{cassirer44} argued that \textit{Erlangen} geometry is also fundamental to vision: we perceive two objectively--different-looking objects as the same because they can be identified by a transformation of space or illumination.
% There are usually nontrivial transformations that leave any given object unchanged; these are called its \emph{symmetry group} and constitute its self-identifications in the ``groupoid of things''.
% Thus, groupoids arise even in social sciences such as psychology; so a language such as \hott\ with groupoids as its basic objects deserves serious philosophical attention.

\section{Conclusion}
\label{sec:conclusion}

There is much more to \hott\ than I have been able to mention in this short chapter, but those aspects I have touched on revolve around a single idea, which generalizes Bishop's set-definition principle: whenever we define a collection of objects, we must also ensure that the identifications and higher identifications between them are correctly specified.
Sometimes the correct identifications arise ``automatically'', such as from the univalence axiom; other times we have to generate new ones, as with higher inductive types.
But in no case must we (or even \emph{can} we) separate those identifications from the objects themselves: with $\infty$-groupoids as basic foundational objects, every collection carries along with itself the appropriate notion of identification between its objects, higher identification between those, and so on.
This can be regarded as the central innovation of \hott, both for mathematics and for philosophy.

% There should be very natural ways to avoid the mistakes of the hole argument using (informal) \hott\ language.

\bibliographystyle{plain}
\bibliography{all}

\end{document}